\documentclass[%
reprint,
superscriptaddress,
amsmath,amssymb,
aps,
]{revtex4-1}

\usepackage{color}
\usepackage{graphicx}
\usepackage{dcolumn}
\usepackage{bm}
\usepackage{amsthm}

\DeclareMathOperator{\arccot}{arccot}

\begin {document}

\title{Aging  arcsine law in Brownian motion and its generalization
}

\author{Takuma Akimoto}
\email{takuma@rs.tus.ac.jp}
\affiliation{%
  Department of Physics, Tokyo University of Science, Noda, Chiba 278-8510, Japan
}%

\author{Toru Sera}

\author{Kosuke Yamato}

\author{Kouji Yano}
\affiliation{%
  Department of Mathematics, Graduate School of Science, Kyoto University, Sakyo-ku, Kyoto 606-8502, Japan
}%



\date{\today}

\begin{abstract}
Classical arcsine law states that fraction of occupation time on the positive or the negative side in Brownian motion 
does not converge to a constant but converges in distribution to the arcsine distribution. Here, we consider 
how a preparation of the system affects the arcsine law, i.e., 
aging of the arcsine law. 
We derive aging distributional theorem for occupation time statistics in Brownian motion, where the 
ratio of time when measurements start to the measurement time plays an important role in determining the shape of the distribution. Furthermore, we show that this result can be generalized as 
aging distributional limit theorem in renewal processes. 
\end{abstract}

\maketitle


Stationarity is one of the most fundamental properties in stochastic processes. 
In equilibrium, physical quantities fluctuate around a constant value and the value is given by the equilibrium ensemble. 
However, statistical properties of physical quantities depend explicitly on time in non-equilibrium processes where 
the characteristic time scale diverges \cite{Bouchaud1992, God2001, Brok2003, Margolin2005, Margolin2006, He2008, Weigel2011, Yamamoto2014, Massignan2014, Miyaguchi2011, Miyaguchi2015, Schulz2013}. In non-stationary stochastic processes, 
aging phenomena are essential, which can be observed by changing the start of the observation time or the total measurement time under the same setup \cite{bouchaud90, Bouchaud1992}. 
In fact, the forward recurrence time 
in renewal processes explicitly depends on 
the time when the observation starts \cite{God2001,Schulz2013}. Furthermore, the mean square displacement (MSD) 
and the diffusion coefficient obtained by
single trajectories depend on the start of the observation as well as the total measurement time in some diffusion processes
\cite{He2008, Weigel2011, Miyaguchi2011, Akimoto2013a, Metzler2014, Yamamoto2014, Massignan2014, Akimoto2014, Miyaguchi2015}.
A typical model that shows aging is a continuous-time random walk (CTRW) with infinite mean waiting time. In the CTRW, 
the MSD increases non-linearly  \cite{metzler00}, i.e., anomalous diffusion, 
\begin{equation}
\langle x(t)^2 \rangle \sim  D_\alpha t^\alpha \quad (t\to \infty), 
\end{equation}
where $x(t)$ is a displacement and $0<\alpha< 1$ characterizes the power-law exponent of the waiting time distribution.  Moreover, it shows 
aging; i.e., the MSD explicitly depends on the start of the observation: 
\begin{equation}
\langle \{x(t_a + t) - x(t_a) \}^2 \rangle \sim  D_\alpha \{ (t_a + t)^\alpha -t_a^\alpha\}
\end{equation}
for $t_a \gg 1$, where $t_a$ is called aging time. 

Aging phenomena are also observed in weakly chaotic dynamical systems such as Pomeau-Manniville map \cite{PM1979, Manneville1980, Akimoto2013b}. 
In weakly chaotic maps, the invariant measure cannot be normalized, i.e., infinite measure \cite{Akimoto2010a}. Moreover,  
the generalized Lyapunov exponent, which characterizes a dynamical instability of the system, 
depends explicitly on the aging time \cite{Akimoto2013b}. In particular, the dynamical instability becomes weak 
when the aging time is increased. 
When the invariant measure of a dynamical system cannot be normalized, the density of a position does not converge to 
the invariant measure. This situation is similar to non-equilibrium processes exhibiting aging.  
In dynamical systems with infinite measures, time-averaged observables do not converge to a constant but converge in distribution in the long-time 
limit \cite{Aaronson1981, Aaronson1997}. In particular, distribution of time averages of $L^1(\mu)$ function, i.e., a function integrable with respect to 
invariant measure $\mu$,  converge to the Mittag-Leffler distribution \cite{Aaronson1981, Aaronson1997}. 
These distributional behaviors of time averages 
are characteristics of infinite ergodic theory, which includes the Mittag-Leffler distribution, the generalized arcsine distribution and another distribution
 \cite{Thaler1998, Thaler2002, TZ2006, Akimoto2008, Akimoto2015,MR3988607,Sera2020}. 

Aging distributional limit theorem in renewal processes, i.e., aging of the Mittag-Leffler distribution, has been studied 
in Refs.~\cite{Schulz2013,Schulz2014} and  is applied to a weakly chaotic 
dynamical system \cite{Akimoto2013b}. However, aging of the arcsine law has not been considered so far. 
In this paper, we consider aging of the arcsine law, which is a distributional theorem of occupation time on the positive side in the Brownian motion. 
We rigorously prove an aging distributional theorem in the Brownian motion. Moreover, we generalize the aging distributional theorem to that 
in renewal processes in the long-time limit. Finally, we demonstrate numerically the aging distributional limit theorem 
using intermittent maps with infinite measures.

{\it Preliminaries.}---Consider 1D Brownian motion starting from the origin. 
This fundamental model of a stochastic process is described by 
 $\dot{x}(t) =  \xi (t)$, 
 where  $\xi (t)$ is a white Gaussian noise:
 $\langle \xi(t) \xi (t') \rangle = \delta (t-t')$. 
 As is well known, the MSD grows as $\langle x(t)^2 \rangle =t$, implying diffusion coefficient $D$ is $D=1/2$.  
  In what follows, we denote the Brownian motion at time $t$ by $B_t$.

Here, we recall the first-passage time (FPT) distribution of a Brownian motion starting from position $x$, the classical arcsine law, and give some notations. 
Let $P_x(s)$ be the probability density function (PDF) of FPT, which is the time when a Brownian motion starting from position $x$ reaches zero for the first time. 
It is known that the PDF  is given by
\begin{equation}
P_x(s)= \frac{x}{s} p(s,x)
\end{equation}
for all $x>0$ and $s>0$ \cite{karatzas2012brownian}, where $p(s,x)$ is the propagator of a Brownian motion, i.e., 
\begin{equation}
 p(s,x) = \frac{1}{\sqrt{2\pi s}} e^{-\frac{x^2}{2s}}
\end{equation}
for $s>0$ and $x\in {\mathbb R}$.

{\it Lemma 1}: For all $t>0$, the distribution of FPT $D_t$,  which is the time when a Brownian motion  reaches zero for the first time after time $t$ passed, 
i.e.,  $D_t \equiv \inf \{ s >0; B_{s+t} =0 \}$ follows $\Pr (D_t > s) = \int_s^\infty \psi_t (u)du$, where 
\begin{equation}
\psi_t (s) = \frac{1}{\pi}  \frac{\sqrt{t}}{\sqrt{s}(s+t)}
\label{frt-brown}
\end{equation}
{\it Proof}: Integrating $P_x(s) p(t,x)$ with respect to $x$, we have
\begin{equation}
\psi_t (s)= \int_{-\infty}^\infty \frac{x}{s} p(s,x) p(t,x) dx = \frac{1}{\pi} \frac{\sqrt{st}}{s(s+t)}. \qed
\end{equation}

\begin{figure}
\includegraphics[width=.9\linewidth, angle=0]{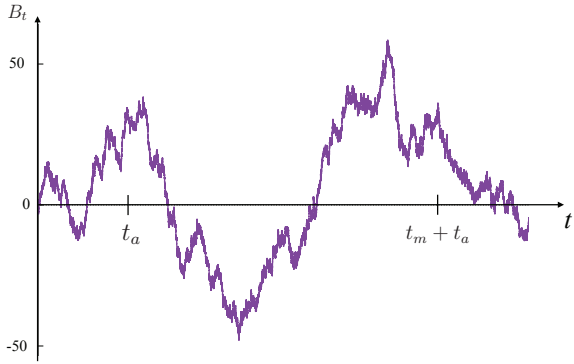}
\caption{Trajectory of Brownian motion with $B_0=0$. In the aging arcsine law, we measure the occupation time from $t_a$ to $t_a+t_m$.}
\label{traj-bm}
\end{figure}

We consider an occupation time $T_+(t)$ that a Brownian motion $B_t$ spends on the positive side until time $t$, i.e., 
\begin{equation}
T_+(t)= \int_0^t 1_{[B_s >0]}ds
\end{equation}
for $t>0$, where 
 $1_{[B_s>0]}=1$ if $B_s>0$ and 0 otherwise. 
The classical arcsine law states that a ratio between an occupation time of a Brownian motion starting from zero 
and measurement time $t_m$  follows the arcsine distribution: 
\begin{equation}
\Pr \left(\frac{T_+(t_m)}{t_m} \leq s\right) = \int_0^s \phi (s')ds' = \frac{2}{\pi} \arcsin \sqrt{s} ,
\end{equation}
where 
\begin{equation}
\phi (s) \equiv \frac{1}{\pi  \sqrt{s (1-s)}} 
\end{equation}
for $0<s<1$. Here, we do not represent the initial position of a Brownian motion explicitly, but it is $B_0=0$. 
By the scaling property of a Brownian motion, this statement is equivalent to the following:
\begin{equation}
\Pr \left(T_+(1) \leq s\right) = \frac{2}{\pi} \arcsin \sqrt{s} .
\end{equation}

{\it Aging arcsine law}.---We introduce the aging time $t_a$, which is a start of measurement [see Fig.~\ref{traj-bm}]. 
Before $t_a$ we do not track the trajectory although the process was 
 started. In other words, a position of a Brownian motion is not the origin when the measurement is started. 

\begin{figure}
\includegraphics[width=.9\linewidth, angle=0]{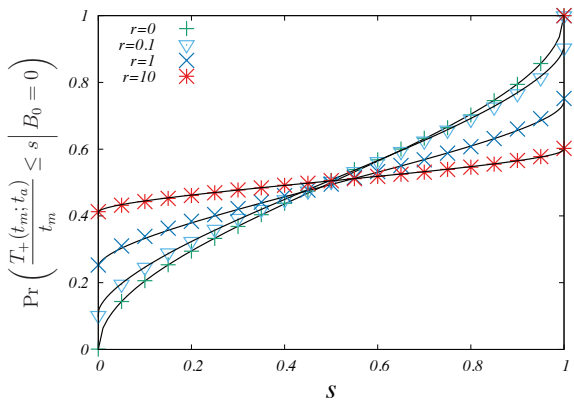}
\caption{Distribution of the ratio of occupation time $T_+(t_m; t_a)$ to $t_m$ in Brownian motion for different 
aging ratio $r$, where measurement time 
$t_m$ is fixed as $t_m=10^3$. Symbols are the results of numerical 
simulations and the solid lines represent our theory, i.e., Eq.~(\ref{prop}).}
\label{ocs-bm}
\end{figure}

\begin{widetext}
{\it Theorem 1}: For all $t_m>0$ and $t_a>0$, the ratio of occupation time $T_+(t_m; t_a) \equiv T_+(t_m+t_a) - T_+(t_a)$ 
to measurement time $t_m$ follows 
 \begin{equation}
 \Pr \left( \left. \frac{T_+(t_m; t_a)}{t_m} \leq s \right| B_0=0  \right) 
 =\int_0^s \phi (r;s')ds' +q(r) + 1_{[s\geq 1]} q(r),
\label{prop}
\end{equation}
where $r\equiv t_a/t_m$ is the aging ratio,  
\begin{equation}
\phi (r;s) = \frac{1}{2\pi^2} \int_0^{1/r} \left\{ \frac{1}{\sqrt{1-s} (1+sv)} + \frac{1}{\sqrt{s} (1+(1-s)v)}
\right\} \frac{dv}{\sqrt{v(1-rv)}} 
\label{pdf-aging-arcsine}
\end{equation}
and
\begin{equation}
q (r) = \frac{1}{2\pi} \int_{1/r}^\infty \frac{dv}{(1+v)\sqrt{v}} 
=\frac{1}{\pi} \arccot (\sqrt{r^{-1}}).
\end{equation}
\end{widetext}
Proof of Theorem 1 is given in the Supplemental Material.  
We note that $\phi(r;s) \to \phi(s)$ for $r\to 0$. In other words, the classical arcsine law is recovered when $t_a \ll t_m$. 
This is consistent with the arcsine law without aging, i.e., $t_a=0$. Figure~\ref{ocs-bm} shows the effect of aging in 
the occupation time statistics.  
In the limit of $r\to 0$, the classical arcsine law is actually recovered. 
Furthermore, 
\begin{equation}
\phi(r;s) \sim c(r) \phi(s)
\end{equation}
for $s\to 0$ and $s\to 1$, where 
\begin{equation}
c(r)= \frac{1}{2\pi} \int_0^{1/r} \frac{dv}{(1+v)\sqrt{v(1-rv)}} .
\end{equation}
Therefore, constant $c(r)$ explicitly depends on  aging ratio $r$. In particular, $c(r)\to 1/2$ and $c(r)\to 0$ for $r\to 0$
and $r\to\infty$, respectively. We note that the classical arcsine law cannot be recovered when limit $s\to 0$ or $s\to 1$ is taken in advance, 
i.e., $c(r)$ does not go to one for $r\to 0$ after $s\to 0$ or $s\to 1$. 
In other words, the limits of $s\to 0$ and $r\to 0$ are not commutative. 

{\it Generalization of the aging arcsine law}.---Here, we generalize our result, i.e., the aging arcsine law, to occupation time statistics in renewal processes \cite{cox, God2001}. 
We consider a two-state process $(R_t)_{t\geq 0}$, 
where the state is described by $+1$ or $-1$ state. Durations for $+1$ and $-1$ states are 
independent and identically distributed random variables. The PDFs of durations for $+1$ and $-1$ states are denoted 
by $\rho_+(\tau)$ and $\rho_-(\tau)$, respectively. We assume that the PDFs follow power-law distributions:
\begin{equation}
\rho_\pm(\tau) \sim A_\pm \tau^{-1-\alpha}\quad (\tau \to \infty). 
\end{equation}
In general, the first duration does not follow $\rho_\pm(\tau)$. However, the following results do not depend 
on the first duration distribution in general.  Therefore, in what follows, we do not specify the initial condition. 
For $0<\alpha < 1$, the mean duration  diverges and the forward recurrence time $D_t^\pm$, which is a time at which 
state changes from $\pm$ to $\mp$, respectively, for the first time after time $t$, shows aging. In particular, the PDF of $D_t^\pm$
depends explicitly on time $t$ \cite{Dynkin1961}. Let us define  $\psi_t^\pm (\tau)$ as
\begin{equation}
\psi_{t}^\pm (\tau) =  p_\pm\frac{\sin \pi \alpha}{\pi} \frac{t^\alpha}{\tau^\alpha (\tau+t)} ,
\end{equation}
where $p_\pm$ is the probability of finding state is $\pm$ at time $t$ and given by $p_\pm =A_\pm/(A_++A_-)$.  
In the limit of $t_m\to\infty$ with being $t_a/t_m=r$ fixed, we have 
\begin{equation}
\Pr\left[ \frac{D_{t_a}^\pm}{t_m} \leq s, R_{t_a} \gtrless 0\right] \to \int_0^s \psi_{r}^\pm (\tau) d\tau. 
\end{equation}
This is consistent with Brownian motion's result, i.e., Eq.~(\ref{frt-brown}), where $\alpha=p_\pm=1/2$ in the Brownian motion. 
 
 \begin{widetext}
In the renewal process, 
the classical arcsine law can be generalized. Occupation time of $+1$ state in the renewal processes follows the generalized arcsine law 
\cite{Lamperti1958, Thaler2002}:
\begin{equation}
\Pr \left( \frac{T_+(t)}{t} \leq s \right) \to \frac{1}{\pi\alpha} \arccot \left( \frac{((1-s)/s)^\alpha}{\beta \sin \pi \alpha}
+\cot \pi \alpha\right)= 
\int_0^s \phi_{\alpha, \beta} (s')ds'\quad (t\to \infty),
\end{equation}
where ${\displaystyle T_+(t) = \int_{0}^{t} 1_{[R_s>0]}ds}$, $\beta=A_-/A_+$ and 
\begin{equation}
\phi_{\alpha, \beta} (s) = \frac{\beta \sin \pi \alpha}{\pi} \frac{s^{\alpha -1}(1-s)^{\alpha -1}}{\beta^2 s^{2\alpha} + 2\beta s^\alpha 
(1-s)^\alpha \cos \pi \alpha + (1-s)^{2\alpha}}.
\end{equation}

{\it Theorem 2}: In the limit of $t_m\to\infty$ with being $t_a/t_m=r$ fixed, 
the ratio of occupation time $T_+(t_m; t_a)$ measured from $t_a$ to $t_m+t_a$ to measurement time $t_m$ follows 
 \begin{equation}
 \Pr \left(  \frac{T_+(t_m; t_a)}{t_m} \leq s  \right) \to \Phi_{\alpha,\beta} (r; s) \equiv \int_0^s \phi_{\alpha, \beta} (r;s')ds' +q_\alpha^-(r) 
+ 1_{[s\geq 1]} q_\alpha^+(r),
\label{tht2}
\end{equation}
where ${\displaystyle T_+(t_m; t_a) = \int_{t_a}^{t_a+t_m} 1_{[R_s>0]}ds}$, 
\begin{equation}
\phi_{\alpha, \beta} (r;s) = \int_0^{s} \psi_{r}^+ (s')  \phi_{\alpha, \beta} \left( \frac{s-s'}{1-s'} \right)  \frac{ds'} {1-s'}+ 
\int_0^{1-s} \psi_{r}^-(s')  \phi_{\alpha, \beta} \left(\frac{s}{1-s'} \right) \frac{ds'}{1-s'} 
\label{pdf-aging-arcsine2}
\end{equation}
and
\begin{equation}
q_\alpha^\pm (r) =  \int_{1}^\infty \psi_r^\pm (\tau)d\tau .
\end{equation}
\end{widetext}
Proof of Theorem 2 is given in the Supplemental Material. 

\begin{figure}
\includegraphics[width=.9\linewidth, angle=0]{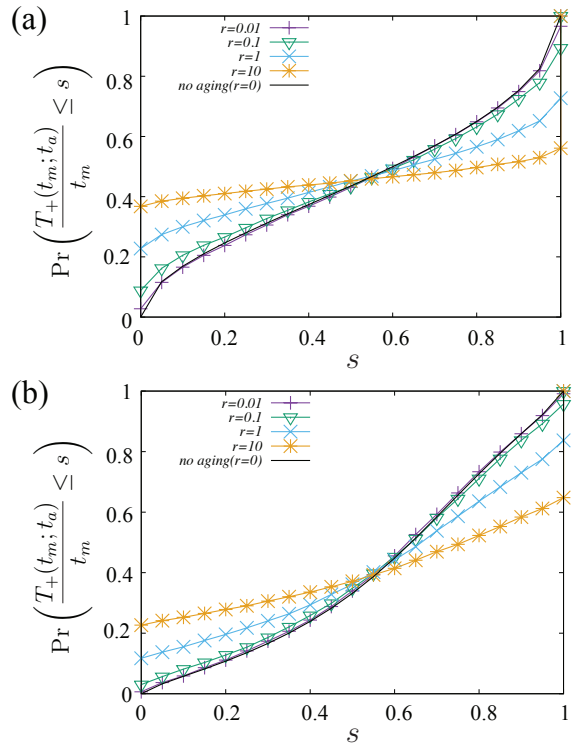}
\caption{Distribution of the ratio of occupation time $T_+(t_m; t_a)$ to measurement time $t_m$  
in the intermittent map, i.e., Eq.~(\ref{smbm}), for different aging ratio $r$, where 
measurement time $t_m$ is fixed as $t_m=10^5$ [(a) $\alpha=0.5$ and (b) $\alpha=0.7$ ($c=0.6$)]. 
Symbols with lines are the results of numerical 
simulations and the solid line represent the generalized arcsine distribution without aging, i.e., $\Phi_{\alpha,\beta}(s) = \int_0^s 
\phi_{\alpha,\beta}(s') ds'$. We used a uniform distribution as the initial distribution.}
\label{ocs-smbm}
\end{figure}

{\it Application of the aging generalized arcsine law to occupation time statistics in intermittent maps}.---Here, we apply the aging distributional limit theorem in renewal processes to occupation time statistics in intermittent maps. One-dimensional 
 map that we consider here is defined on $[0,1]$, i.e., $T(x): [0,1] \to [0,1]$:
 \begin{equation}
 T(x) = \left\{
 \begin{array}{ll}
 x + (1-c) \left(\dfrac{x}{c} \right)^{1+1/\alpha} \quad &x\in [0,c]\\
 \\
 x - c \left(\dfrac{1-x}{1-c} \right)^{1+1/\alpha} \quad &x \in (c,1], 
 \end{array}
\right.
\label{smbm}
 \end{equation}
where $c$ ($0<c<1$) is a parameter characterizing a skewness of the map 
and $0<\alpha <1$ \cite{Akimoto2008}. There are two indifferent fixed points 
at $x=0$ and 1, i.e., $T(0)=0$ and $T(1)=0$ with $T'(0)=T'(1)=1$. With the aid of the chaotic behaviors outside the two indifferent 
fixed points, durations on $[0,c]$ or $(c,1]$ are considered to be independent and identically distributed random variables. 
Moreover, the duration distributions follow  a power-law \cite{Akimoto2008, Akimoto2015}. Therefore, the aging distributional limit
theorem can be applied to the occupation time statistics in the intermittent map. In the case of no aging, the ordinary generalized 
arcsine law is shown \cite{Thaler2002}, where parameter $\beta$ is given by
\begin{equation}
\beta = \frac{\alpha + c}{\alpha + 1 -c} \left( \frac{1-c}{c}\right)^{2\alpha}.
\end{equation}
Figure~\ref{ocs-smbm} shows the distribution 
of the ratio of occupation time $T_+(t_m; t_a)$ to $t_m$ on $[0,c]$. 
The shape of the distribution strongly depends on aging ratio $r$. Moreover, 
the generalized arcsine distribution can be recovered for small $r$. This is because the generalized arcsine distribution 
is obtained by substituting $r=0$ and $\psi_r^+ (s) =p_+ \delta (s)$ and $\psi_r^- (s) =p_- \delta (s)$ in Eq.~(\ref{pdf-aging-arcsine2}).

{\it Conclusion}.---We have shown aging distributional theorem of occupation time on the positive side in Brownian motion. 
FPT distribution $P_x(s)$ of Brownian motion starting from $x$ is a key distribution to derive the theorem. 
The distribution of the occupation time is described by aging ratio $r$. 
The classical arcsine law is recovered when aging time $t_a$ is much smaller than measurement time $t_m$, 
i.e., $r\to 0$. We have also shown that the aging arcsine law is generalized to the occupation time distribution in renewal processes
under the limits of $t_a \to \infty$ and $t_m \to \infty$. The ordinary generalized arcsine law is also recovered in the limit 
of $r\to 0$. Finally, this generalized aging arcsine law can be successfully applied to the occupation time statistics 
in  intermittent maps with infinite invariant measures.

%

\clearpage

\begin{widetext}
\appendix

{\bf Supplemental Material for ``Aging  arcsine law in Brownian motion and its generalization"}

Here, we give proofs of Theorem 1 and 2.

\section{Proof of Theorem 1}
{\it Proof}: By the scaling property of the Brownian motion, statistical properties of $B_{st}$ are the same as those of 
$\sqrt{t}B_s$. It follows that statistical properties of occupation time $T_+(t_m; t_a)/t_m$ are the same 
as those of $T_+(r+1) - T_+(r)$ because
\begin{equation}
\frac{T_+(t_m; t_a)}{t_m} = \frac{1}{t_m} \int_{t_a}^{t_a+t_m} 1_{[B_s>0]} ds = \int_{r}^{r+1} 1_{[B_{st_m}>0]} ds.
\end{equation}
First, we consider case $B_{t_a}>0$. Using the scaling property, we have
\begin{eqnarray}
\Pr \left( \frac{T_+(t_m; t_a)}{t_m} \leq s , B_{t_a}>0 \right) = \Pr \left(T_+(r+1) - T_+(r) \leq s, B_r>0 \right) .
\end{eqnarray}
Since the probability of $B_{t_a}>0$ is 1/2, 
\begin{eqnarray}
\Pr \left( \frac{T_+(t_m; t_a)}{t_m} \leq s , B_{t_a}>0 \right) 
= \frac{1}{2} \int_0^s \psi_r(s') \Pr (T_+ (1-s') \leq s-s')   ds' 
\end{eqnarray}
for $s<1$ and 
\begin{eqnarray}
\Pr \left( \frac{T_+(t_m; t_a)}{t_m} =1 , B_{t_a}>0 \right) 
= \int_1^\infty \frac{\psi_r (s')}{2}ds' .
\end{eqnarray}
By a change of variables, we obtain 
\begin{eqnarray}
\Pr \left( \frac{T_+(t_m; t_a)}{t_m} \leq s , B_{t_a}>0 \right) 
&=& \frac{1}{2} \int_0^s \psi_r(s') \Pr \left(T_+ (1) \leq \frac{s-s'}{1-s'} \right)   ds' + 1_{[s\geq 1]} q(r) \nonumber\\
&=& \frac{1}{2\pi^2} \int_0^{\frac{s}{r}} \frac{dy}{\sqrt{y}(1+y)} \int_0^{\frac{s-ry}{1-ry}} \frac{du}{\sqrt{u(1-u)}}  
+ 1_{[s\geq 1]} q(r) \nonumber\\
&=& \int_0^{s}  \frac{1}{2\pi^2} \frac{dv}{\sqrt{1-v}} \int_0^{\frac{1}{r}} \frac{du}{\sqrt{u (1-ru)}(1+vu)} + 1_{[s\geq 1]} q(r).
\end{eqnarray}
By a similar calculation, we have
\begin{eqnarray}
\Pr \left( \frac{T_+(t_m; t_a)}{t_m} \geq s , B_{t_a}<0 \right) 
= \int_s^{1}  \frac{1}{2\pi^2} \frac{dv}{\sqrt{v}} \int_0^{\frac{1}{r}} \frac{du}{\sqrt{u (1-ru)}(1+(1-v)u)}
\end{eqnarray}
 for $s>0$. 
For $B_{t_a}<0$ and $s=0$, the probability is 
\begin{eqnarray}
\Pr \left( \frac{T_+(t_m; t_a)}{t_m} = 0 , B_{t_a}<0 \right) = \frac{1}{2} \int_1^\infty \psi_r (s')ds'. 
\end{eqnarray}
It follows that aging arcsine distribution is given by Eq.~(\ref{prop}) and the PDF $\phi (r;s)$ is given by Eq.~(\ref{pdf-aging-arcsine}).
\qed

\section{Proof of Theorem 2}

{\it Proof}: By a scaling 
argument, aging occupation time statistics can be obtained by a similar way in the Brownian motion. By a change of variables, 
we have
\begin{equation}
\frac{T_+(t_m; t_a)}{t_m} 
= \tilde{T}_+(r+1) - \tilde{T}_+(r) ,
\end{equation}
where ${\displaystyle \tilde{T}_+(r) = \int_{0}^{r} 1_{[R_{st_m}>0]} ds}$. 
We note that limits $t_a \gg 1$ and $t_m\gg 1$ are necessary to derive the distribution of occupation time 
in renewal processes, which is different from the arcsine law in the Brownian motion. 
For $R_{t_a}>0$ and $t_m\gg 1$ and $t_a=r t_m \gg 1$, we have 
\begin{eqnarray}
\Pr \left( \frac{T_+(t_m; t_a)}{t_m} \leq s , R_{t_a}>0 \right) &=& \Pr \left(\tilde{T}_+(r+1) - \tilde{T}_+(r) \leq s, R_{t_a}>0 \right) .
\end{eqnarray}
By a similar calculation as in the aging arcsine law, we obtain 
\begin{eqnarray}
\Pr \left( \frac{T_+(t_m; t_a)}{t_m} \leq s , R_{t_a}>0 \right)
\to \int_0^{s} dv \int_0^{v} \frac{\psi_{r}^+ (s')}{1-s'}  \phi_{\alpha, \beta} \left( \frac{v-s'}{1-s'} \right)  ds' + 1_{[s\geq 1]} q_\alpha^+(r).
\end{eqnarray}
Similarly, 
\begin{eqnarray}
\Pr \left( \frac{T_+(t_m; t_a)}{t_m} \geq s , R_{t_a}<0 \right) 
\to  \int_{s}^{1} dv \int_0^{1-v} \psi_{r}^-(s')  \phi_{\alpha, \beta} \left(\frac{v}{1-s'} \right) \frac{1}{1-s'} ds' 
\end{eqnarray}
for $s>0$ and 
\begin{eqnarray}
\Pr \left( \frac{T_+(t_m; t_a)}{t_m} = 0 , R_{t_m}<0 \right) \to \int_1^\infty \psi_r^- (s')ds'. 
\end{eqnarray}
It follows that aging arcsine distribution is given by Eq.~(\ref{tht2}) and the PDF $\phi_{\alpha,\beta}  (r;s)$ is given by Eq.~(\ref{pdf-aging-arcsine2}).\qed

\end{widetext}

\end{document}